# W Sequences and the Distribution of Primes in Short Interval

Shaohua Zhang

School of Mathematics in Shandong University, Jinan, Shandong, People's Republic of China 250100

E-mail address: shaohuazhang@mail.sdu.edu.cn

**Abstract:** Based on Euclid's algorithm, we find a kind of special sequences which play an interesting role in the study of primes. We call them W Sequences. They not only ties up the distribution of primes in short interval but also enables us to give new weakened forms of many classical problems in Number Theory. The object of this paper is to provide a brief introduction and preliminary analysis on this kind of special sequences.

## I Introduction

As we know, using Euclid's algorithm, one can solve the following two problems:

**Problem 1:** Let $a_1,...,a_m$ and $b$ be any positive integers. Find an algorithm to determine whether $b$ can be represented by $a_1,...,a_m$.

**Problem 2:** Let $a_1,...,a_m$ be any positive integers. Find an algorithm to determine whether there is an integer $a_i$ among $a_1,...,a_m$ such that $a_i$ is relatively prime with all of the others.

Generalize Problem **1** to the case over the unique factorization domain $F[x_1,...,x_n]$, where $F$ is a field. Let $f_1,...,f_m$ and $g$ be polynomials in $F[x_1,...,x_n]$. Find an algorithm to determine whether $g$ belongs to the ideal generated by $f_1,...,f_m$. As we know, this interesting generalization leads to the invention of Gröbner bases.

Problem **2** leads to the invention of a kind of special sequences which play an interesting role in the study of primes and enable us to give new weakened forms of many classical problems which are open in Number Theory.

As we know, the probability that $k$ random integers are pairwise relatively prime is $1/\zeta(k)$, where $\zeta(.)$ is Riemann's Zeta function. The probability that there is an integer $a_i$ among several random integers $a_1,...,a_m$ such that $a_i$ is relatively prime with all of the others is $(1/\zeta(2))^{m-1}$. For such an integer sequence $a_1,...,a_m$, in which there is an integer $a_i$ such that $a_i$ is relatively prime with all of the others, we are very interested in it and would like to give its proprietary definition for the convenience of expression. For the details, see Section 2.

## 2 W Sequences and their properties

**Definition of W sequence:** For any integer $n > 1$, the sequence of integers $0 < a_1 < ... < a_n$ is called a W sequence, if $\exists r (1 \leq r \leq n)$ such that $a_r$ and each of the rest numbers are coprime.

And here, $a_r$ is called a W number of the sequence.

Let's consider firstly W sequences in the case of consecutive positive integers. The following Lemma 1 and theorems are simple, we omit the proofs.

**Lemma 1:** For $n > 1$, $m \geq 1$, if $m+1, m+2, ..., m+n$ is a W sequence, let $m+r$ be a W number and $p$ be the smallest prime factor of $m+r$, where $1 \leq r \leq n$, then $p \geq \frac{n+1}{2}$.

**Theorem 1:** If there is a prime in $\{m+1, m+2, ..., m+n\}$, then this sequence is a W sequence.

**Theorem 2:** For $m \geq 1$, that there exists a prime in the interval $(m^2, (m+1)^2)$ is equivalent with that $m^2+1, ..., m^2+2m$ is a W sequence.

**Theorem 3:** For $\varepsilon > 0$, let $m$ be a sufficiently large natural number, then that there is a prime in the interval $(m, m+m^{\frac{1}{2}+\varepsilon}]$ is equivalent with that $m+1, ..., m+[m^{\frac{1}{2}+\varepsilon}]$ is a W sequence.

By the theorem 1 above, in order to determine whether a consecutive positive integer sequence is a W sequence, it is enough to consider the case of consecutive composite numbers. On consecutive composite numbers, Grimm [2] made an important conjecture in 1969.

Denote the largest integer $n$ in $m+1, m+2, ..., m+n$ by $h(m)$ such that no one of $m+1, m+2, ..., m+n$ is relatively prime with all of the others. Cramér's conjecture [1] and the result of Pillai [3] imply $17 \leq h(m) << (\log m)^2$ by the theorem 1.

From the aforementioned results, one can see that the W sequences tie up the distribution of primes in short interval. Unfortunately, it is not easy to prove that a sequence is a W sequence. As for W sequences in the case of non-consecutive positive integers, they are more intractability. However, they enable us to get new necessary conditions of many classical problems. As an example, if for $m > 1, n > 1$, there is a prime in the Fibonacci sequence $f_{m+1}, f_{m+2}, ..., f_{m+n}$, then $f_{m+1}, f_{m+2}, ..., f_{m+n}$ is a W sequence.

For another example, for positive integers $a$ and $b$ with $(a,b)=1$ and $2<b$, if there is a prime in $a, a+b, a+2b, ..., a+(b-1)b$, then this sequence is a W sequence. Thus, in order to determine whether the sequence $a, a+b, a+2b, ..., a+(b-1)b$ is a W sequence, it is enough to consider the case that $a, a+b, a+2b, ..., a+(b-1)b$ are all composite numbers. It leads to the generalization of Grimm's conjecture. In [5], we give a possible generalization.

In 2008, by studying W sequences in the case of non-consecutive positive integers, we give a new weakened form of Goldbach's Conjecture and reveal the internal relationship between Goldbach's conjecture and the least prime in an arithmetic progression [4].

Finally, we hope that people can find a generic criterion for determining whether a sequence is a W sequence or not in the near future. Let's ask several interesting questions to close this paper.

Question 1: For any $m>3$, is $1^2+1, 2^2+1, ..., m^2+1$ a W sequence?

Question 2: For any integers $a$ and $b$ with $(a,b)=1$, $2<b$ and $0<a<b$, is the sequence $a, a+b, a+2b, ..., a+(b-1)b$ a W sequence?

Question 3: Let $p_i$ be the $i$th odd prime and let $m>1$. Is $2+p_1, 2+p_2, ..., 2+p_m$ a W sequence?

## Acknowledgements

Thank my advisor Professor Xiaoyun Wang for her valuable help. Thank the Institute for Advanced Study in Tsinghua University for providing me with excellent conditions. This work was partially supported by the National Basic Research Program (973) of China (No. 2007CB807902) and the Natural Science Foundation of Shandong Province (No.Y2008G23).